\documentclass[papera4,12pt]{article}
\usepackage[english]{babel}
\usepackage[letterpaper,top=2cm,bottom=2cm,left=3cm,right=3cm,marginparwidth=1.75cm]{geometry}
\usepackage{amsmath}
\usepackage{graphicx}
\usepackage[english]{babel}
\usepackage[colorlinks=true, allcolors=blue]{hyperref }
\usepackage{amsthm}
\usepackage{amsfonts}
\usepackage{hyperref}
\newtheorem{thm}{Theorem}[section]

\newtheorem{lem}[thm]{Lemma}

\numberwithin{equation}{section}
\begin{document}
\begin{center}
{\bf {\Large Pell and Pell-Lucas numbers as sums of two Jacobsthal numbers }}
\vspace{8mm}

{\Large \bf  Ahmed Gaber}\\
\vspace{3mm}
  Department of Mathematics, University of Ain Shams \\
Faculty of science, Egypt \\
e-mails:  \href{mailto:a.gaber@sci.asu.edu.eg}{\url{a.gaber@sci.asu.edu.eg}}
\end{center}
\vspace{10mm}

\begin{abstract}
 We solve the two Diophantine equations $P_k=J_n+J_m$ and $Q_k=J_n+J_m$
where $\left\lbrace P_{k}\right\rbrace_{k\geq0}$, $\left\lbrace
Q_{k}\right\rbrace_{k\geq0}$ and $\left\lbrace J_{k}\right\rbrace_{k\geq0}$ are the sequences of Pell numbers, Pell-Lucas numbers and Jacobsthal numbers, respectively.  The main tool is the theory of   linear forms in logarithms.
\end{abstract}
{\bf Keywords:} Pell sequence, Pell-Lucas sequence, Jacobsthal sequence, Linear forms in logarithms.\\
{\bf 2020 Mathematics Subject Classification:} 11B39,11D72,11J70.
\vspace{5mm}
\section{Introduction}
The Pell sequence $\lbrace{P_k\rbrace}$ is defined recursively by $P_0=0, P_1=1$ and $P_{n+1}=2P_{n}+P_{n-1}$ for $n\geq 1$ . A few terms of this sequence are\\
$$0, 1, 2, 5, 12, 29, 70, 169, 408, 985,...$$
  Pell-Lucas numbers are defined by  $Q_0=2$, $Q_1=2$, and $Q_n+1=2Q_{n}+Q_{n-1}$ for $n\geq 1$. Its first few terms are
$$2, 2, 6, 14, 34, 82, 198, 478, 1154, 2786,...$$
For some recent works related to Diophantine equations which includes Pell and Pell-Lucas  numbers, see \cite{2}, \cite{1} and \cite{3}.\\
Jacobsthal sequence is defined by $J_0=0$, $J_1=1$, and $J_n+1=J_{n}+2J_{n-1}$ for $n\geq 1$. Its first few terms are
$$0, 1, 1, 3, 5, 11, 21, 43, 85, 171,...$$

The aim of this paper is to study the two  diophantine equations
  \begin{equation}\label{eqn1}
    P_k=J_n+J_m,
\end{equation}
and
  \begin{equation}\label{eqn100}
    Q_k=J_n+J_m.
\end{equation}
 The complete sets of solutions are provided in the following theorems.
\begin{thm}\label{thm1}
The only non-negative solutions $(k,n,m)$  which satisfy  Eq.(\ref{eqn1}) and $n\geq m$ are
\begin{center}
  $(1,1,0), (2,2,1), (2,2,2), (3,4,0), (4,5,1), (4,5,2).$
\end{center}

\end{thm}
\begin{thm}\label{thm40}
The only non-negative solutions $(k,n,m)$  which satisfy  Eq.(\ref{eqn100}) and $n\geq m$ are
\begin{equation*}
  (0,1,1),(1,1,1),(1,2,1),(1,2,2),(2,3,3),(2,4,1),(2,4,2),(3,5,3),(0,2,1),(1,2,1).
\end{equation*}
\end{thm}
\section{Preliminary results}

\subsection{Pell and  Pell-Lucas sequences }
 The characteristic equation of Pell and  Pell-Lucas Sequences is
\begin{equation*}
    \Psi(x):=x^2-2x-1=0,
   \end{equation*}
and the Binet formula of the Pell and  Pell-Lucas Sequences are, assuming  $\gamma=1+\sqrt{2}$ and  $\delta=1-\sqrt{2}$
\begin{equation}
P_k=\frac{\gamma^{k}-\delta^{k}}{2\sqrt{2}} \hspace{0.35 cm} \text{for all}\hspace{0.25 cm}k\geq 0.
\end{equation}

 and\\
 \begin{equation}
Q_k=\gamma^{k}+\delta^{k} \hspace{0.35 cm} \text{for all}\hspace{0.25 cm}k\geq0.
\end{equation}
Also, by induction one can prove that
\begin{equation}\label{eqn8}
    \gamma^{k-2} \leq P_k \leq \gamma^{k-1} \hspace{0.25cm}\text{holds for all}\hspace{0.25cm} k\geq 1.
\end{equation}
and
\begin{equation}\label{eqn201}
    \gamma^{k-1}< Q_{k} < \gamma^{k+1} \hspace{0.25cm}\text{holds for all}\hspace{0.25cm} k\geq 2.
\end{equation}

\subsection{Jacobsthal sequence}
  The characteristic equation of the Jacobsthal sequence is
\begin{equation*}
    P(x)=x^2-x-2.
\end{equation*}
Its Binet formula is
\begin{equation}
    J_n=\frac{2^n-(-1)^n}{3}.
\end{equation}
It is known that
\begin{equation}\label{eqn10}
    2^{n-2} \leq J_n \leq 2^{n-1}\hspace{0.25cm} \text{for all }\hspace{0.25cm} n \geq 1.
\end{equation}
Basic properties of Jacobsthal numbers can be found in \cite{6} and \cite{7}.
\subsection{Linear forms in logarithms.  }
Let $\eta $ be an algebraic number of degree $d$ with minimal polynomial over $\mathbb{Z}$
\begin{equation*}
    a_{0}x^{d}+a_{1}x^{d-1}+...+a_{d}=a_{0}\prod_{i=1}^{d} (x-\eta^{(i)}),
\end{equation*}
where the leading coefficient $a_{0}$ is positive and the $\eta^{(i)}$'s are the conjugates of $\eta$. Then the logarithmic height of $\eta$ is given by
\begin{equation*}
    h(\eta):=\frac{1}{d}\left(\log a_{0} +\sum_{i=1}^{d} \log \left( \max \left\{\left| \eta^{(i)}\right|,1\right\}\right)\right)
\end{equation*}
The following are some of the properties of the logarithmic height function $h(.)$:
\begin{eqnarray*}
h\left(\eta_{1}\pm \eta_{2}\right)&\leq& h\left(\eta_{1}\right)+h\left(\eta_{2}\right)+\log 2;\\
h\left(\eta_{1}\eta_{2}^{\pm 1}\right)&\leq& h\left(\eta_{1}\right)+h\left(\eta_{1}\right);\\
h\left(\eta^{s}\right)&=&\left|s\right|h\left(\eta\right) \hspace{0.2cm} \left(s\in \mathbb{Z}\right).
\end{eqnarray*}
The proof of the following result can be found in \cite{10}.
\begin{thm}\label{thm3}
Let $\eta _{1},...,\eta _{l}$ be positive real algebraic numbers in a real algebraic number field $\mathbb{L}\subset\mathbb{R}$ of degree $D$, $b_{1},...,b_{l}$ be a non zero integers, and assume that
\begin{equation*}
    \Lambda_{1}:=\eta_{1}^{b_{1}}...\eta_{l}^{b_{l}}-1\neq 0.
\end{equation*}
Then,
\begin{equation}
    \log \left| \Lambda\right|>-1.4\cdot 30^{l+3}\cdot l^{4.5}\cdot D^{2}\cdot(1+\log D)\cdot (1+\log B) A_1...A_{l},
\end{equation}
where
\begin{equation*}
B\geq \max \lbrace{\left|b_{1}\right|,...,\left|b_{l}\right|\rbrace},
\end{equation*}
and
\begin{equation*}
    A_{i}\geq \max \lbrace{Dh\left(\eta_{i}\right),\left|\log \eta_{i}\right|,0.16\rbrace}, \hspace{0.2cm} \text{for all} \hspace{0.3cm} i=1,...,l.
\end{equation*}
\end{thm}
\subsection{Dujella and Peth{\"o} reduction lemma}
 Let $X$ be a real number. Set $\left||X\right||:=min\lbrace{\left|X-n\right|:n\in \mathbb{Z}\rbrace}$. In \cite{5} Dujella and Peth{\"o} proved the following important reduction result .
\begin{lem}\label{lem2.4}
Let $M$ be a positive integer. Let $\tau,\mu,A >0, B>1$ be given real numbers. Assume $\frac{p}{q}$ is a convergent of $\tau$ such that $q >6M$ and $\epsilon : = ||\mu q||-M||\tau q||>0.$ If $(n,m,\omega) $ is a positive solution to the inequality
\begin{equation*}
    0<\left| n\tau -m+\mu\right|<\frac{A}{B^{\omega}}
\end{equation*}
with $n \leq M,$ then
\begin{equation*}
    \omega < \frac{\log \left(\frac{Aq}{\epsilon}\right)}{\log B}.
\end{equation*}
\end{lem}
\subsection{Legendre theorem}
The following theorem is due to  Legendre and will be used in some cases of our investigation of Pell-Lucas numbers that are expressible as the sum of two Jacobsthal numbers.  Further details can be found in \cite{456}.
\begin{thm}\label{thm20}
Let $x$ be a real number, let $p,q \in \mathbb{Z}$ and let $x=[a_{0},a_{1},...]$.
If
\begin{equation*}
    \left|\frac{p}{q}-x\right|<\frac{1}{2 q^2}
\end{equation*}
then $\frac{p}{q}$ is a convergent continued fraction of $x$. Furthermore, let M and n be a non negative integers such that $q_n > M$. Put $b=\max\lbrace{a_i:i=0,1,2,..,n\rbrace}$ then,
\begin{equation*}
    \frac{1}{\left(b+2\right)q^2}<\left|\frac{p}{q}-x\right|
\end{equation*}
\end{thm}
\section{Proof of Theorem \ref{thm1}}

\subsection{Bounding $m$, $n$ and $k$.}
Applying the inequalities (\ref{eqn10}) and (\ref{eqn8}) to establish the relationship between $k$ and $n$  Then we get,
\begin{equation}
    \gamma^{k-3} \leq P_k \leq 2^{n} \hspace{0.2cm}  \text{and}\hspace{0.2cm}  2^{n-2} \leq P_k \leq \gamma^{k-1}.
\end{equation}
These implies that,
\begin{equation}
    (n-2)\frac{\log 2}{\log \gamma}+1 \leq k \leq n \frac{\log 2}{\log \gamma}+2.
\end{equation}
We can consider $k < 2n$. Using the Binet formulas of the Pell and Jacobsthal sequences in (\ref{eqn1}). We get
 \begin{equation}\label{eqn200}
  \frac{\gamma^{k}-\delta^{k}}{2\sqrt{2}}   =\frac{2^n-(-1)^n}{3}+\frac{2^m-(-1)^m}{3}.
\end{equation}
Then,
\begin{equation}
    \left| \frac{\gamma^{k}}{2\sqrt{2}}-\frac{2^{n}}{3}\right| = \left| \frac{2^{m}}{3}-\frac{\left((-1)^{n}+(-1)^{m}\right)}{3}+\frac{\delta^{k}}{2\sqrt{2}}\right|.
\end{equation}
This implies that
\begin{equation}
    \left| \frac{\gamma^{k}}{2\sqrt{2}}-\frac{2^{n}}{3}\right| < \frac{4\cdot 2^m}{3}.
\end{equation}
Thus,
\begin{equation}\label{eqn16}
   \left| \frac{3\gamma^{k}2^{-n}}{2\sqrt{2} }-1\right|<\frac{4}{2^{n-m}}.
\end{equation}
Let
\begin{equation*}
    \Lambda_{1}=\frac{3\gamma^{k}2^{-n}}{2\sqrt{2} }-1,\hspace{0.2cm} l=3,\hspace{0.2cm} \eta_{1}= \frac{3}{2\sqrt{2}},\hspace{0.2cm} \eta_{2}=\gamma,\hspace{0.2cm} \eta_{3}=2,\hspace{0.2cm} b_{1}=1,\hspace{0.2cm} b_{2}=k,\hspace{0.2cm} b_{3}=-n.
\end{equation*}
  If $\Lambda_{1}= 0,$ then  $3 \gamma^{k}=2^{n}.2\sqrt{2}.$ Consider the automorphism $\sigma$ such that $\sigma( \gamma)=\delta.$ Then $\left|3\delta^{k}\right|=2^n.2\sqrt{2}$. But $\left|3\delta^{k}\right|<3$, then $2^n.2\sqrt{2}<3$ which is a contradiction. So, $\Lambda_{1}\neq 0.$

Take $\mathbb{L}=\mathbb{Q}(\gamma).$ Then $D=2.$ The logarithmic heights are
\begin{equation*}
    h(\eta_{1})\leq h(3)+h(2\sqrt{2}) \leq \log 3 +\frac{3}{2} \log 2;
\end{equation*}
\begin{equation*}
    h(\eta_{2})=\frac{1}{2} \log \gamma;
\end{equation*}
\begin{equation*}
    h(\eta_{3})=\log 2.
\end{equation*}
Taking,
\begin{equation*}
A_1=2\log 3 +3\log 2, \hspace{0.2cm} A_2=\log \gamma, \hspace{0.2cm} \text{and} \hspace{0.2cm} A_3=2\log 2, B=2n
    \end{equation*}
    and applying Matveev’s Theorem(\ref{thm3}), we obtain
    \begin{equation*}
      \log \left| \Lambda_{1}\right|> -1.4\times 30^{6}\times 3^{4.5}\times 4\times (1+\log 2)(1+\log 2n)(2\log 3 +3\log 2)(2\log 2 \log \gamma),
    \end{equation*}
    then,
    \begin{equation}\label{eqn17}
      \log \left| \Lambda_{1}\right|>-6\times 10^{12}(1+\log 2n).   \end{equation}
      Also, from (\ref{eqn16}) we have,
      \begin{equation}\label{eqn18}
         \log \left| \Lambda_{1}\right|<\log 4+(m-n)\log 2.
      \end{equation}
      Thus by comparing inequalities in (\ref{eqn17}) and (\ref{eqn18}), we get
      \begin{equation}\label{eqn19}
      (n-m)\log 2 -\log 4 < 6\times 10^{12}(1+\log 2n).
 \end{equation}
 Hence,
 \begin{equation}\label{eqn20}
 m\log 2 >n\log 2- 6\times 10^{12}(1+\log 2n)-\log 4.
 \end{equation}
  Eq.(\ref{eqn200}) can be written as,
 \begin{equation}
     \frac{\gamma^{k}}{2\sqrt{2} }-\frac{2^{n}(1+2^{m-n})}{3}=\frac{\delta^{k}}{2\sqrt{2} }-\frac{(-1)^{n}-(-1)^{m}}{3},
 \end{equation}
 Therefore,
 \begin{equation}
     \left|\frac{3\gamma^{k}2^{-n}}{2\sqrt{2}(1+2^{m-n})}-1\right|=\left|\frac{3\cdot2^{-n}}{2\sqrt{2}(1+2^{m-n})}\left(\frac{\delta^{k}}{2\sqrt{2} }-\frac{(-1)^{n}-(-1)^{m}}{3}\right)\right|.
 \end{equation}
 Hence,
 \begin{equation}\label{eqn23}
     \left|\frac{3\gamma^{k}2^{-n}}{2\sqrt{2}(1+2^{m-n})}-1\right|<\frac{5}{2^m}.
 \end{equation}
 Let $\Lambda_{2}=\frac{3\gamma^{k}2^{-n}}{2\sqrt{2}(1+2^{m-n})}-1.$
 Then,
 \begin{equation}\label{eqn24}
     \log\left|\Lambda_{2}\right|< \log 5 -m\log 2.
 \end{equation}
 Let
\begin{equation*}
     \eta_{1}=\frac{3}{2\sqrt{2}(1+2^{m-n})},\hspace{0.2cm} \eta_{2}=\gamma,\hspace{0.2cm} \eta_{3}=2,\hspace{0.2cm} l=3,\hspace{0.2cm} b_{1}=1,\hspace{0.2cm} b_{2}=k,\hspace{0.2cm} b_{3}=-n, B=2n.
\end{equation*}
First we show that $\Lambda_{2}\neq 0.$ If $\Lambda_{2}= 0,$ then  $3 \gamma^{k}=2\sqrt{2}(2^{n}+2^{m}).$ Consider the automorphism $\sigma$ such that $\sigma(a\gamma)=\delta.$ Then $\left|3\delta^{k}\right|=2\sqrt{2}(2^{n}+2^{m}).$ But $\left|3\delta^{k}\right|<3$,   which is a contradiction.
Then we take $\mathbb{L}=\mathbb{Q}(\gamma),$ for which $D=2.$ We compute the logarithmic heights as follows:
\begin{equation*}
    h(\eta_{1})\leq h(3)+h(2)+h(2\sqrt{2})+h(1+2^{m-n}) \leq \log 3 + \frac{5}{2} \log 2+(n-m)\log 2 ;
    \end{equation*}
\begin{equation*}
    h(\eta_{2})=\frac{1}{2} \log \gamma;
\end{equation*}
\begin{equation*}
    h(\eta_{3})=\log 2.
\end{equation*}
We take,
\begin{equation*}
A_1=2\log 3 +2(n-m)\log 2+5\log 2, \hspace{0.2cm} A_2=\log \gamma, \hspace{0.2cm} \text{and} \hspace{0.2cm} A_3=2\log 2.
    \end{equation*}
  Then, by sMatveev’s Theorem,  we get
 \begin{equation*}
      \log \left| \Lambda_{2}\right|> -1.4\times 30^{6}\times 3^{4.5}\times 4\times (1+\log 2)(1+\log 2n)(2\log 3 +2(n-m)\log 2+5\log 2)(2\log 2\log \gamma),
    \end{equation*}
 Using Eqs.(\ref{eqn19}),(\ref{eqn20}),(\ref{eqn24}) and direct computations,  we find
 \begin{equation}\label{eqn25}
  n \log 2 < 6\times 10^{13}(1+\log 2n)+24\times 10^{27}(1+\log2n)^{2}+3.
   \end{equation}
  We deduce that
\begin{equation}\label{eqn26}
   n<2\times 10^{29}.
 \end{equation}

\subsection{Reducing bound on $n$}

Now, we use the reduction lemma to  reduce the upper bound on $n.$
     \begin{equation*}
         \Gamma_{1}=\log(\frac{3}{2\sqrt{2}})+k \log \gamma-n\log 2.
     \end{equation*}
  Eq.(\ref{eqn16}) gives
  \begin{equation}\label{eqn28}
      \Lambda_{1}=e^{\Gamma_{1}}-1<\frac{4}{2^{m-n}}<\frac{1}{4},
  \end{equation}
 which implies that
 \begin{equation}\label{eqn29}
     \left|\Gamma_{1}\right|<\frac{1}{2}.
 \end{equation}
 Then $\left|\Gamma_{1}\right|<2\left|e^{\Gamma_{1}}-1\right|.$ Therefore we get
 \begin{equation}
     \left|\Gamma_{1}\right|<\frac{8}{2^{n-m}}.
 \end{equation}
 We observe that $\Gamma_{1}\neq0$ since $\Lambda_{1}\neq 0.$ Then
 \begin{equation}\label{eqn31}
     0 < \left|\frac{\log(\frac{3}{2\sqrt{2}})}{\log 2}-n+k\left(\frac{\log \gamma }{\log 2}\right)\right| < \frac{8}{2^{n-m} \log 2} <\frac{12}{2^{n-m}}.
 \end{equation}
 We apply lemma (\ref{lem2.4}) with $M=4\times 10^{29}$ ($M>2n>k$) , $\tau = \frac{\log \alpha }{\log 2} $ , $\mu=\frac{\log(\frac{3}{2\sqrt{2}})}{\log 2} $, $A=12$, $B=2$. Write $\tau$ as a continued fraction $[a_{0},a_{1},...]$ we get\\ $q_{65}=2427228558134035529638808203392547 >  6M$. We compute
 \begin{equation*}
     \epsilon = \left|| \mu q_{65}\right||-M\left||\tau q_{65}\right||>0.1.
 \end{equation*}
 Thus by lemma (\ref{lem2.4}), we get  $n-m <118$.
 Now we put
 \begin{equation*}
         \Gamma_{2}=\log \left(\frac{3}{2\sqrt{2}(1+2^{m-n})}\right)+k \log \gamma-n\log 2.
     \end{equation*}
  Then we have from equation(\ref{eqn23}) that, for $m\geq 5$,
  \begin{equation}\label{eqn36}
      \Lambda_{2}=e^{\Gamma_{2}}-1<\frac{5}{2^{m}}<\frac{1}{4},
  \end{equation}
 which implies that
 \begin{equation}\label{eqn37}
     \left|\Gamma_{2}\right|<\frac{1}{2}.
 \end{equation}
 Then $\left|\Gamma_{2}\right|<2\left|e^{\Gamma_{2}}-1\right|.$ Therefore we get
 \begin{equation}
     \left|\Gamma_{2}\right|<\frac{10}{2^m}.
 \end{equation}
 We observe that $\Gamma_{2}\neq0$ since $\Lambda_{2}\neq 0.$ Then
 \begin{equation}\label{eqn31}
     0 < \left|\frac{\log \left(\frac{3}{2\sqrt{2}(1+2^{m-n})}\right)}{\log 2}-n+k\left(\frac{\log \gamma }{\log 2}\right)\right| <\frac{15}{2^{m}}.
 \end{equation}
 We apply lemma (\ref{lem2.4}) with $M=4\times 10^{29}$ ($M>2n>k$),  $\tau = \frac{\log \gamma }{\log 2} $ , $\mu=\frac{\log (\frac{3}{2\sqrt{2}(1+2^{m-n})})}{\log 2} $ , $A=15$ and  $B=2$. It can be seen that $q_{65}=2427228558134035529638808203392547 >  6M$. Computing every $\epsilon$ such that $n-m<118$, we get
 \begin{equation*}
     \epsilon = \left|| \mu q_{65}\right||-M\left||\tau q_{65}\right||>0.01.
 \end{equation*}
 Thus by lemma (\ref{lem2.4}), we get  $m < 122$. So, $n<240$ and $k<480$ . Solving Eq.(\ref{eqn1}) for $m < 122$,  $n<240$ and $k<480$, we got the indicated solutions in Theorem (\ref{thm1}). The proof is complete.

\section{Proof of Theorem \ref{thm40}}

By symmetry of Eq.(\ref{eqn100}), we assume that $n\geq m$.
\subsection{bounding $m$, $n$ and $k$.}
By (\ref{eqn201}) and (\ref{eqn10}),  we have
\begin{equation}
    \alpha^{k-2} \leq R_k \leq 2^{n} \hspace{0.2cm}  \text{and}\hspace{0.2cm}  2^{n-2} \leq R_k \leq \alpha^{k+1}.
\end{equation}
These imply that,
\begin{equation}
    (n-2)\frac{\log 2}{\log \gamma}-1 \leq k \leq (n-1) \frac{\log 2}{\log \alpha}+1.
\end{equation}
We take $k < 2n.$ Replace the Pell-Lucas and Jacobsthal sequences in (\ref{eqn100}) by their  Binet formulas as follows:
\begin{equation}\label{eqn500}
   \gamma^{k}+\delta^{k}=\frac{2^n-(-1)^n}{3}+\frac{2^m-(-1)^m}{3}.
\end{equation}
Then,
\begin{equation}
    \left| \gamma^{k}-\frac{2^{n}}{3}\right| = \left| \frac{2^{m}}{3}-\frac{\left((-1)^{n}+(-1)^{m}\right)}{3}-\delta^{k}\right|.
\end{equation}
Therefore,
\begin{equation*}
    \left| \gamma^{k}-\frac{2^{n}}{3}\right| \leq \frac{2^{m}}{3}+\frac{2}{3}+ \left|\delta^k\right|
\end{equation*}
Then
\begin{equation}\label{eqn42}
     \left| \gamma^{k}-\frac{2^{n}}{3}\right|<\frac{4\cdot 2^m}{3}
\end{equation}
Thus,
\begin{equation}\label{eqn43}
   \left| \frac{3\gamma^{k}}{2^{n}}-1\right|<\frac{4}{2^{n-m}}.
\end{equation}
Consider the following:
\begin{equation*}
    \Lambda_{3}=3 \gamma^{k} 2^{-n}-1,\hspace{0.2cm} l=3,\hspace{0.2cm} \eta_{1}=3,\hspace{0.2cm} \eta_{2}=\gamma,\hspace{0.2cm} \eta_{3}=2,\hspace{0.2cm} b_{1}=1,\hspace{0.2cm} b_{2}=k,\hspace{0.2cm} b_{3}=-n.
\end{equation*}
We show that $\Lambda_{3}\neq 0.$ If $\Lambda_{3}= 0,$ then  $3 \gamma^{k}=2^{n}.$ Consider the automorphism $\sigma$ such that $\sigma(\gamma)=\delta.$ Then $\left|3\delta^{k}\right|=2^n.$ But $\left|3\delta^{k}\right|<3$, then $2^n<3$ which is a contradiction.\\
Take $\mathbb{L}=\mathbb{Q}(\gamma),$ for which $D=2.$ Then,
\begin{equation*}
    h(\eta_{1})=\log 3 ;
\end{equation*}
\begin{equation*}
    h(\eta_{2})=\frac{1}{2} \log \gamma;
\end{equation*}
\begin{equation*}
    h(\eta_{3})=\log 2.
\end{equation*}
We take,
\begin{equation*}
A_1=2\log 3 , \hspace{0.2cm} A_2=\log \gamma, \hspace{0.2cm} \text{and} \hspace{0.2cm} A_3=2\log 2.
    \end{equation*}
   Let $B=2n.$ Then Theorem(\ref{thm3}) shows that
    \begin{equation*}
      \log \left| \Lambda_{3}\right|> -1.4\times 30^{6}\times 3^{4.5}\times 4\times (1+\log 2)(1+\log 2n)(2\log 3)(2\log 2\log \gamma).
    \end{equation*}
    Consequently,
    \begin{equation}\label{eqn44}
      \log \left| \Lambda_{3}\right|>-3\times 10^{12}(1+\log 2n).
      \end{equation}
      Then also from (\ref{eqn43}) we have,
      \begin{equation}\label{eqn45}
         \log \left| \Lambda_{3}\right|<\log 4+ (m-n)\log 2.
      \end{equation}
      Thus by comparing inequalities in (\ref{eqn44}) and (\ref{eqn45}) we get,
      \begin{equation}\label{eqn46}
      (n-m)\log 2 -\log 6 < 3\times 10^{12}(1+\log 2n).
 \end{equation}
 Hence,
 \begin{equation}\label{eqn47}
 m\log 2 >n\log 2- 3\times 10^{12}(1+\log 2n)-\log 4.
 \end{equation}
 Eq.(\ref{eqn500}) is equivalent to
 \begin{equation}\label{eqn48}
     \gamma^{k}-\frac{2^{n}(1+2^{m-n})}{3}=-\delta^{k}-\frac{(-1)^{n}-(-1)^{m}}{3}.
 \end{equation}
 So,
 \begin{equation}\label{eqn49}
     \left|\frac{3\gamma^{k}2^{-n}}{1+2^{m-n}}-1\right|=\left|\frac{3\cdot2^{-n}}{1+2^{m-n}}\left(-\delta^{k}-\frac{(-1)^{n}-(-1)^{m}}{3}\right)\right|.
 \end{equation}
 Then,
 \begin{equation}\label{eqn50}
     \left|\frac{3\gamma^{k}2^{-n}}{1+2^{m-n}}-1\right|<\frac{5}{2^m}.
 \end{equation}
 Let $\Lambda_{4}=\frac{3}{1+2^{m-n}} \gamma^{k} 2^{-n}-1.$
 Hence,
 \begin{equation}\label{eqn51}
     \log\left|\Lambda_{4}\right|< \log 5 -m\log 2.
 \end{equation}
 Set
\begin{equation*}
     \eta_{1}=\frac{3}{1+2^{m-n}},\hspace{0.2cm} \eta_{2}=\gamma,\hspace{0.2cm} \eta_{3}=2,\hspace{0.2cm} l=3,\hspace{0.2cm} b_{1}=1,\hspace{0.2cm} b_{2}=k,\hspace{0.2cm} b_{3}=-n.
\end{equation*}
First we show that $\Lambda_{4}\neq 0.$ If $\Lambda_{4}= 0,$ then  $3 \gamma^{k}=2^{n}+2^{m}.$ Consider the automorphism $\sigma$ such that $\sigma(\gamma)=\delta.$ Then $\left|3\delta^{k}\right|=2^n+2^{m}.$ But $\left|3\delta^{k}\right|<3$, then $2^n+2^{m}<3$ which is a contradictio.\\ 
Let $\mathbb{L}=\mathbb{Q}(\gamma),$ for which $D=2.$ Then,
\begin{equation*}
    h(\eta_{1}) \leq \log 3 +(n-m)\log 2 +\log 2;
    \end{equation*}
\begin{equation*}
    h(\eta_{2})=\frac{1}{2} \log \gamma;
\end{equation*}
\begin{equation*}
    h(\eta_{3})=\log 2.
\end{equation*}
We take,
\begin{equation*}
A_1=2\log 3 +2(n-m)\log 2+2\log 2, \hspace{0.2cm} A_2=\log \gamma, \hspace{0.2cm} \text{and} \hspace{0.2cm} A_3=2\log 2 and B=2n
    \end{equation*}
      Then
    \begin{equation*}
      \log \left| \Lambda_{4}\right|> -1.4\times 30^{6}\times 3^{4.5}\times 4\times (1+\log 2)(1+\log 2n)(2\log 3+2\log 2+2(n-m)\log 2)(2\log 2\log \gamma),
    \end{equation*}
    Then using equations (\ref{eqn46}), (\ref{eqn47}) and (\ref{eqn51}) with some computation we get,
    \begin{equation}\label{eqn52}
      n \log 2 < 11\times 10^{12}(1+\log 2n)+12\times 10^{24}(1+\log2n)^{2}+3.
      \end{equation}
     Hence,
     \begin{equation}\label{eqn53}
         n<3\times 10^{28}.
     \end{equation}

     \subsection{Reducing bound on $n$.}
      Assume that $n-m \geq 5$ and let
     \begin{equation*}
         \Gamma_{3}=\log(3)+k \log \gamma-n\log 2.
     \end{equation*}
  By Eq.(\ref{eqn43}), we have
  \begin{equation}\label{eqn54}
      \Lambda_{3}=e^{\Gamma_{3}}-1<\frac{4}{2^{m-n}}<\frac{1}{4},
  \end{equation}
 So,
 \begin{equation}\label{eqn55}
     \left|\Gamma_{3}\right|<\frac{1}{2}.
 \end{equation}
 Then, $\left|\Gamma_{3}\right|<2\left|e^{\Gamma_{3}}-1\right|.$ Therefore, we have
 \begin{equation}\label{eqn56}
     \left|\Gamma_{3}\right|<\frac{8}{2^{n-m}}.
 \end{equation}
 We observe that $\Gamma_{3}\neq0$ since $\Lambda_{3}\neq 0.$ Then
 \begin{equation}\label{eqn57}
     0 < \left|\frac{\log 3}{\log 2}-n+k\left(\frac{\log \gamma }{\log 2}\right)\right| < \frac{8}{2^{n-m} \log 2} <\frac{12}{2^{n-m}}.
 \end{equation}
 We apply lemma (\ref{lem2.4}) with $M=6\times 10^{28}$ ($M>2n>k$) , $\tau = \frac{\log \gamma }{\log 2} $ , $\mu=\frac{\log 3}{\log 2} $, $A=12$, $B=2$. Considering the continued fraction of  $\tau$,  we find that $q_{65}>6M$. We compute
 \begin{equation*}
     \epsilon = \left|| \mu q_{65}\right||-M\left||\tau q_{65}\right||>0.3.
 \end{equation*}
 Thus, by lemma (\ref{lem2.4}), we get  $n-m < 117$.
 Set
 \begin{equation*}
         \Gamma_{4}=\log \left(\frac{3}{1+2^{m-n}}\right)+k \log \gamma-n\log 2.
     \end{equation*}
    and let $m>5$. Then we have from equation(\ref{eqn50}) that
  \begin{equation}\label{eqn58}
      \Lambda_{4}=e^{\Gamma_{4}}-1<\frac{5}{2^{m}}<\frac{1}{4}.
  \end{equation}
 We conclude that
 \begin{equation}\label{eqn59}
     \left|\Gamma_{4}\right|<\frac{1}{2}.
 \end{equation}
 Thus $\left|\Gamma_{4}\right|<2\left|e^{\Gamma_{4}}-1\right|.$ Therefore we get
 \begin{equation}\label{eqn60}
     \left|\Gamma_{4}\right|<\frac{10}{2^m}.
 \end{equation}
 We observe that $\Gamma_{4}\neq0$. So
 \begin{equation}\label{eqn61}
     0 < \left|\frac{\log \left(\frac{3}{1+2^{m-n}}\right)}{\log 2}-n+k\left(\frac{\log \gamma }{\log 2}\right)\right| <\frac{15}{2^{m}}.
 \end{equation}
 We apply lemma (\ref{lem2.4}) with $M=6\times 10^{28}$ ($M>2n>k$),  $\tau = \frac{\log \gamma }{\log 2} $ , $\mu=\frac{\log \left(\frac{3}{1+2^{m-n}}\right)}{\log 2} $ , $A=15$, $B=2$. We have\\ $q_{65}>6M$. We consider the values of $\epsilon$ in two cases\\
 $\mathbf{Case ~I}$: if $n-m<117$ and $n-m\neq 1$

 \begin{equation*}
     \epsilon = \left|| \mu q_{65}\right||-M\left||\tau q_{65}\right||>0.01.
 \end{equation*}
 Thus by lemma (\ref{lem2.4}), we get  $m < 122$ so $n<239$ and $k<478$ . \\

 $\mathbf{Case ~II}$: $n-m=1$ we get $\epsilon$ always negative. So we solve equation(\ref{eqn100}) if $n-m=1$. In this case  equation(\ref{eqn100}) can be written as $Q_k=J_{m}+J_{m+1}$ and can be reduced as
 \begin{equation}\label{eqn62}
     Q_{k}=2^{m}.
 \end{equation}
 Then $k<2m$ and  from Eq.(\ref{eqn53}) we get $m<3\times 10^{28}$. As before, We can prove that
 \begin{equation*}
     \alpha^{k}2^{-m}-1<\frac{1}{2^{m}},
 \end{equation*}
 This gives, for $m\geq 3$, that
 \begin{equation*}
     \left|k\frac{\log \gamma}{\log 2}-m\right|<\frac{4}{2^m}< \frac{1}{4}.
 \end{equation*}
 Using the relation $16m < 2^{m}$ for $m\geq 7$, we deduce that $\frac{4}{2^m}<\frac{1}{2k^2}$. Then $\left|\frac{\log \gamma}{\log 2}-\frac{m}{k}\right|<\frac{1}{2k^{2}}$.
So, by Legendre's theorem,  $\frac{m}{k}$ is a convergent of $\frac{\log \gamma}{\log 2}$. Using $k<M$ and some computations we find that
 \begin{equation*}
     q_{53}< M < q_{54} \hspace{0.3cm} \text{and} \hspace{0.3cm} b:=\max \{a_{i}: i=0,1,2,...,54\}<2.
 \end{equation*}
 Consequently,
 \begin{equation*}
     \frac{1}{\left(2+2\right) k}< \frac{4}{2^m}.
 \end{equation*}
 Thus
 \begin{equation*}
     2^{m}<16\cdot 6\cdot 10^{28}.
 \end{equation*}
 Then $m\leq 100.$ Solutions of Eq.(\ref{eqn62}) for $m<100$  and Eq.(\ref{eqn100})  for $m < 122$, $n<239$ and $k<478$  completes the proof of Theorem (\ref{thm40}).\\

\makeatletter
\renewcommand{\@biblabel}[1]{[#1]\hfill}
\makeatother

\end{document}